\documentclass{amsproc}

\newtheorem{theorem}{Theorem}[section]
\newtheorem{lemma}[theorem]{Lemma}
\newtheorem{proposition}[theorem]{Proposition}

\theoremstyle{definition}
\newtheorem{definition}[theorem]{Definition}

\theoremstyle{remark}
\newtheorem{remark}[theorem]{Remark}

\numberwithin{equation}{section}

\begin{document}

\newcommand{\ea}{\mbox{{\bf a}}}

\newcommand{\eu}{\mbox{{\bf u}}}

\newcommand{\ueu}{\underline{\eu}}

\newcommand{\ueo}{\overline{u}}

\newcommand{\oeu}{\overline{\eu}}

\newcommand{\ew}{\mbox{{\bf w}}}

\newcommand{\ef}{\mbox{{\bf f}}}

\newcommand{\eF}{\mbox{{\bf F}}}

\newcommand{\eC}{\mbox{{\bf C}}}

\newcommand{\en}{\mbox{{\bf n}}}

\newcommand{\eT}{\mbox{{\bf T}}}

\newcommand{\eL}{\mbox{{\bf L}}}

\newcommand{\eR}{\mbox{{\bf R}}}

\newcommand{\eV}{\mbox{{\bf V}}}

\newcommand{\eU}{\mbox{{\bf U}}}

\newcommand{\ev}{\mbox{{\bf v}}}

\newcommand{\eve}{\mbox{{\bf e}}}

\newcommand{\uev}{\underline{\ev}}

\newcommand{\eY}{\mbox{{\bf Y}}}

\newcommand{\eK}{\mbox{{\bf K}}}

\newcommand{\eP}{\mbox{{\bf P}}}

\newcommand{\eS}{\mbox{{\bf S}}}

\newcommand{\eJ}{\mbox{{\bf J}}}

\newcommand{\eB}{\mbox{{\bf B}}}

\newcommand{\eH}{\mbox{{\bf H}}}

\newcommand{\leb}{\mathcal{ L}^{n}}

\newcommand{\eI}{\mathcal{ I}}

\newcommand{\eE}{\mathcal{ E}}

\newcommand{\hen}{\mathcal{H}^{n-1}}

\newcommand{\eBV}{\mbox{{\bf BV}}}

\newcommand{\eA}{\mbox{{\bf A}}}

\newcommand{\eSBV}{\mbox{{\bf SBV}}}

\newcommand{\eBD}{\mbox{{\bf BD}}}

\newcommand{\eSBD}{\mbox{{\bf SBD}}}

\newcommand{\ecs}{\mbox{{\bf X}}}

\newcommand{\eg}{\mbox{{\bf g}}}

\newcommand{\paromega}{\partial \Omega}

\newcommand{\gau}{\Gamma_{u}}

\newcommand{\gaf}{\Gamma_{f}}

\newcommand{\sig}{{\bf \sigma}}

\newcommand{\gac}{\Gamma_{\mbox{{\bf c}}}}

\newcommand{\deu}{\dot{\eu}}

\newcommand{\dueu}{\underline{\deu}}

\newcommand{\dev}{\dot{\ev}}

\newcommand{\duev}{\underline{\dev}}

\newcommand{\weak}{\stackrel{w}{\approx}}

\newcommand{\mild}{\stackrel{m}{\approx}}

\newcommand{\strong}{\stackrel{s}{\approx}}

\newcommand{\weakdown}{\rightharpoondown}

\newcommand{\opg}{\stackrel{\mathfrak{g}}{\cdot}}

\newcommand{\opunu}{\stackrel{1}{\cdot}}
\newcommand{\opdoi}{\stackrel{2}{\cdot}}

\newcommand{\opn}{\stackrel{\mathfrak{n}}{\cdot}}
\newcommand{\opx}{\stackrel{x}{\cdot}}

\newcommand{\tr}{\ \mbox{tr}}

\newcommand{\Ad}{\ \mbox{Ad}}

\newcommand{\ad}{\ \mbox{ad}}

\title{Self-similar dilatation structures and automata}

\author{Marius Buliga}
\address{"Simion Stoilow" 
Institute of Mathematics, Romanian Academy, 
P.O. BOX 1-764, RO 014700, 
Bucure\c sti, Romania}
\email{Marius.Buliga@imar.ro}


\date{14.09.2007}

\begin{abstract}
We show that on the boundary of the dyadic tree, any self-similar
 dilatation structure induces a web of interacting  automata.
\end{abstract}

\maketitle


\section*{Introduction}

In this paper we continue the study of dilatation structures, introduced in  
\cite{buligadil1}.

A dilatation structure $(X,d,\delta)$  describes  the approximate 
self-similarity of the 
metric space $(X,d)$.  Metric spaces which admit  strong dilatation structures 
(definition 
\ref{defweakstrong}) have metric tangent spaces at any point (theorem 
7 \cite{buligadil1}). By theorems 8, 10 \cite{buligadil1}, 
 any such metric tangent space has an algebraic structure of a conical group. 
 Particular examples of conical  groups are Carnot groups, that is  simply connected Lie 
groups  whose Lie algebra admits a positive graduation.

Here we are concerned with dilatation structures on ultrametric spaces. 
The special case considered is the boundary of the infinite dyadic tree,  
topologically the same as the middle-thirds Cantor set. This is also 
the space of infinite words over the alphabet $X = \left\{ 0, 1\right\}$.  
Self-similar dilatation structures are introduced  and 
studied on this space. 

We show that on the boundary of the dyadic tree, any self-similar 
dilatation structure is described by a web of interacting  automata. This 
is achieved in theorems \ref{tstruc} and \ref{th2}. These theorems are 
analytical in nature, but they admit an easy interpretation in terms 
of automata by using classical results as  theorem \ref{tcont} and 
proposition \ref{psync}. Due to the limitations in length of the paper, 
we leave this straightforward interpretation, as well as examples, for 
a further paper (but see also the  slow-paced  introduction into the subject 
\cite{buligadil2}).

The subject is relevant for applications to the hot topic of 
self-similar groups of isometries of the dyadic tree (for an introduction 
into self-similar groups see \cite{bargrinek}).

\section{Words and the Cantor middle-thirds set}

Let $X$ be a finite, non empty  set. The elements of $X$ are called letters.  
The collection of words of finite length in the alphabet $X$ is denoted by $X^{*}$. The empty word 
$\emptyset$ is an element of $X^{*}$. 

 The length of any word $w\in X^{*}$,  
$\displaystyle w = a_{1} ... a_{m}$,  $\displaystyle a_{k} \in X$ 
for all $k=1, ... , m $,  
is denoted by $\mid w \mid = m$. 
The set of words which are infinite at right is denoted by 
$$X^{\omega} = \left\{ f \  \mid \ \   f: \mathbb{N}^{*}\rightarrow X \right\} = X^{\mathbb{N}^{*}} \quad . $$
Concatenation of words is naturally defined. If $\displaystyle q_{1},q_{2} \in X^{*}$ and 
$w \in X^{\omega}$ then $\displaystyle q_{1}q_{2} \in X^{*}$ and $q_{1}w \in X^{\omega}$. 

The shift map $\displaystyle s : X^{\omega} \rightarrow X^{\omega}$ is defined 
by $\displaystyle  w = w_{1} \, s(w) $, 
for any word $\displaystyle w \in X^{\omega}$. For any $\displaystyle k \in \mathbb{N}^{*}$ we define  
$\displaystyle [w]_{k} \in X^{k} \subset X^{*}$, $\displaystyle \left\{ w \right\}_{k} \in X^{\omega}$ by 
$$ w =  [w]_{k} \,Ês^{k}(w)  \quad , \quad \left\{ w \right\}_{k}  = s^{k}(w) \quad . $$

The topology on $\displaystyle X^{\omega}$ is generated by cylindrical sets $qX^{\omega}$, 
for all $q\in X^{*}$.  The topological space $\displaystyle X^{\omega}$  is compact. 

To any $q\in X^{*}$ is associated a continuous injective transformation 
$\hat{q}:X^{\omega}\rightarrow X^{\omega}$, $\hat{q}(w) = qw$. The semigroup 
$X^{*}$ (with respect to concatenation) can be identified with the semigroup 
(with respect to function composition) of these transformations. This 
semigroup is obviously generated by $X$. The empty word $\emptyset$  
corresponds to the identity function. 

The dyadic tree $\mathcal{T}$  is the infinite rooted planar binary tree.  
Any node has two descendants. The nodes are  coded by elements of 
$\displaystyle X^{*}$, $X = \left\{ 0,1\right\}$. The root is coded by the 
empty word and if a node is coded by $x\in X^{*}$ then its left hand side 
descendant has the code $x0$ and its 
right hand side descendant has the code $x1$. We shall therefore identify the 
dyadic tree with $\displaystyle X^{*}$ and we put on the dyadic tree the 
natural (ultrametric) distance on $\displaystyle X^{*}$. The boundary (or 
the set of ends) of the dyadic tree is then the same as the compact 
ultrametric space $\displaystyle X^{\omega}$.

\section{Automata}
In this section we use the same notations as \cite{grigorchuk3}. 

\begin{definition}
An  (asynchronous) automaton is an oriented set $\displaystyle (X_{I}, X_{O}, 
Q, \pi, \lambda)$, with: 
\begin{enumerate}
\item[(a)] $\displaystyle X_{I}, X_{O}$ are finite sets, called the input and 
output alphabets, 
\item[(b)] $Q$ is a set of internal states of the automaton, 
\item[(c)] $\pi$ is the transition function, $\displaystyle \pi: X_{I} \times 
Q \rightarrow Q$, 
\item[(d)] $\lambda$ is the output function, $\displaystyle 
\lambda: X_{I} \times Q \rightarrow X_{O}^{*}$. 
\end{enumerate}
If $\lambda$ takes values in $\displaystyle X_{O}$ then the automaton is 
called synchronous. 
\end{definition}

The functions $\lambda$ and $\pi$ can be continued to the set $\displaystyle 
X_{I}^{*} \times Q$ by: $\pi(\emptyset, q)\, = \, q$, $\lambda(\emptyset, q) \, 
 = \, \emptyset$, 
$$\pi(xw, q) \, = \, \pi(w, \pi(x,q)) \quad , \quad \lambda(xw, q) \, = 
\,  \lambda(x,q) \lambda(w, \pi(x,q))$$ 
for any $\displaystyle x \in X_{I}, q \in Q$ and any $\displaystyle 
w \in X_{I}^{*}$. 

An automaton is nondegenerate if the functions $\lambda$ and $\pi$ can be 
uniquely extended by the previous formul\ae to 
$\displaystyle X_{I}^{\omega} \times Q$.  

To any nondegenerated automaton $\displaystyle (X_{I}, X_{O}, 
Q, \pi, \lambda)$ and any $q \in Q$ is associated the function 
$\displaystyle \lambda(\cdot , q) : X_{I}^{\omega} \rightarrow  
X_{O}^{\omega}$. The following is theorem 2.4 \cite{grigorchuk3}. 

\begin{theorem}
The mapping $\displaystyle f:  X_{I}^{\omega} \rightarrow  
X_{O}^{\omega}$ is continuous if and only if it is defined by a certain 
nondegenerate asynchronous automaton. 
\label{tcont}
\end{theorem}

The proof given in \cite{grigorchuk3} is interesting to read because it 
provides a construction of an automaton which defines the continuous function
$f$.

\section{Isometries of  the dyadic tree}

An isomorphism of $\mathcal{T}$ is just an invertible transformation which preserves the structure of the 
tree. It is well known that isometries  of $\displaystyle (X^{\omega}, d)$ are the same as isometries of 
$\mathcal{T}$.  

Let $\displaystyle A \in Isom(X^{\omega}, d)$ be such an isometry. For any finite word $\displaystyle 
q\in X^{*}$ we may define $\displaystyle A_{q} \in Isom(X^{\omega}, d)$ by 
$$A(qw) = A(q) \, A_{q}(w)$$
for any $\displaystyle w \in X^{\omega}$. Note that in the previous relation $A(q)$ makes sense because 
$A$ is also an isometry of $\mathcal{T}$. 

The following description of isometries of the dyadic tree in terms of automata 
can be deduced from an equivalent formulation of proposition 3.1
\cite{grigorchuk3} (see also proposition 2.18 \cite{grigorchuk3}). 

\begin{proposition}
A function $\displaystyle X^{\omega} \rightarrow X^{\omega}$ is an isometry 
of the dyadic tree if and only if it is generated by a synchronous automaton 
with $\displaystyle X_{I} = X_{O} = X$. 
\label{psync}
\end{proposition}

\section{Motivation: linear structure in terms of dilatations}

 For the normed, real, finite dimensional vector space $\mathbb{V}$, the dilatation based 
at $x$, of coefficient $\varepsilon>0$, is the function 
$$\delta^{x}_{\varepsilon}: \mathbb{V} \rightarrow \mathbb{V} \quad , \quad 
\delta^{x}_{\varepsilon} y = x + \varepsilon (-x+y) \quad . $$
For fixed $x$ the dilatations based at $x$ form a one parameter group which
 contracts any bounded neighbourhood of $x$ to a point, uniformly with respect 
to $x$.

The algebraic structure of $\displaystyle \mathbb{V}$ is encoded in 
dilatations. Indeed,  using dilatations we can recover the operation of addition and multiplication by scalars. 
 
For $\displaystyle x,u,v \in \mathbb{V}$ and $\varepsilon>0$ 
define 
$$\Delta_{\varepsilon}^{x}(u,v) = \delta_{\varepsilon^{-1}}^{\delta_{\varepsilon}^{x} u}
 \delta^{x}_{\varepsilon} v \quad , \quad 
\Sigma_{\varepsilon}^{x}(u,v) = \delta_{\varepsilon^{-1}}^{x} \delta_{\varepsilon}^{\delta_{\varepsilon}^{x} u}
 (v) \quad , \quad inv^{x}_{\varepsilon}(u) =  \delta_{\varepsilon^{-1}}^{\delta_{\varepsilon}^{x} u} x \quad . $$

The meaning of this functions becomes clear if we take the limit as 
 $\varepsilon \rightarrow 0$ of these expressions: 
\begin{equation}
\lim_{\varepsilon\rightarrow 0} \Delta_{\varepsilon}^{x}(u,v) = 
\Delta^{x}(u,v) = x+(-u+v) \quad 
\label{exad}
\end{equation}
$$\lim_{\varepsilon\rightarrow 0} \Sigma_{\varepsilon}^{x}(u,v) = \Sigma^{x}(u,v) = u+(-x+v) \quad ,$$
$$\lim_{\varepsilon\rightarrow 0} inv^{x}_{\varepsilon}(u) = inv^{x}(u) = x-u+x \quad , $$
uniform with respect to $x,u,v$ in bounded sets. 
The function $\displaystyle  \Sigma^{x}(\cdot,\cdot)$ is a group operation, namely the addition operation 
translated such that the neutral element is $x$. Thus, for $x=0$, we recover the group operation. The function 
$\displaystyle inv^{x}(\cdot)$ is the inverse function, and $\displaystyle  \Delta^{x}(\cdot,\cdot)$ is the difference function.

Dilatations behave well with respect to the  distance $d$ induced by the norm, in the following sense: 
for any $\displaystyle x,u,v\in \mathbb{V}$ and any $\varepsilon>0$ we have 
\begin{equation}
\frac{1}{\varepsilon} d( \delta_{\varepsilon}^{x} u , 
\delta^{x}_{\varepsilon} v ) = d(u,v) 
\label{coned}
\end{equation}
This shows that from the metric point of view the space $\displaystyle 
(\mathbb{V}, d)$ is a metric cone, that is $\displaystyle 
(\mathbb{V}, d)$ looks the same at all scales. 

Affine continuous  transformations $A:\mathbb{V} \rightarrow \mathbb{V}$ admit the following 
description in terms of dilatations. (We could dispense of continuity hypothesis in this situation, but 
we want to illustrate a general point of view, described further in the paper).  

\begin{proposition}
A continuous transformation  $A:\mathbb{V} \rightarrow \mathbb{V}$ is affine if and only if for any 
$\varepsilon \in  (0,1)$, $x,y \in \mathbb{V}$ we have 
\begin{equation}
 A \delta_{\varepsilon}^{x} y \ = \ \delta_{\varepsilon}^{Ax} Ay \quad . 
 \label{eq1proplin} 
 \end{equation}
\label{1proplin}
\end{proposition}

The proof is a straightforward consequence of representation formul{\ae}               for the addition, difference and 
inverse operations  in terms of dilatations.

Further on we shall take the dilatations as basic data associated to an
ultrametric space. In order to understand our aim we describe it as follows: 
we shall study a particular ultrametric space (the infinite dyadic tree) 
as if we study the vector space $\mathbb{V}$  by using 
 only the distance $d$ and the dilatations $\displaystyle 
\delta^{x}_{\varepsilon}$ for all $x \in X$ and $\varepsilon > 0$. 

We shall call a triple  $(X,d, \delta)$  
a dilatation structure (see further definition \ref{defweakstrong}), where 
$(X,d)$ is a locally compact metric space  and $\delta$ is a    collection of
 dilatations of the metric space $(X,d)$. 
We shall  add
some compatibility relations between the distance $d$ and dilatations 
$\delta$, which will prescribe: 
\begin{enumerate}
\item[-] the behaviour of the distance 
with respect to dilatations, for example some   form of relation 
(\ref{coned}), 
\item[-] the interaction between dilatations, for example the existence of 
the limit from the left hand side of relation (\ref{exad}). 
\end{enumerate}

\section{Dilatation structures}

This section contains the axioms of a dilatation structure, introduced in 
Buliga \cite{buligadil1}.

\subsection{Notations}

Let $\Gamma$ be  a topological separated commutative group  endowed with a continuous group morphism 
$$\nu : \Gamma \rightarrow (0,+\infty)$$ with $\displaystyle \inf \nu(\Gamma)  =  0$. Here $(0,+\infty)$ is 
taken as a group with multiplication. The neutral element of $\Gamma$ is denoted by $1$. We use the multiplicative notation for the operation in $\Gamma$. 

The morphism $\nu$ defines an invariant topological filter on $\Gamma$ (equivalently, an end). Indeed, 
this is the filter generated by the open sets $\displaystyle \nu^{-1}(0,a)$, $a>0$. From now on 
we shall name this topological filter (end) by "0" and we shall write $\varepsilon \in \Gamma \rightarrow 
0$ for $\nu(\varepsilon)\in (0,+\infty) \rightarrow 0$. 

The set $\displaystyle \Gamma_{1} = \nu^{-1}(0,1] $ is a semigroup. We note $\displaystyle 
\bar{\Gamma}_{1}= \Gamma_{1} \cup \left\{ 0\right\}$
On the set $\displaystyle 
\bar{\Gamma}= \Gamma \cup \left\{ 0\right\}$ we extend the operation on $\Gamma$ by adding the rules  
$00=0$ and $\varepsilon 0 = 0$ for any $\varepsilon \in \Gamma$. This is in agreement with the invariance 
of the end $0$ with respect to translations in $\Gamma$.

We shall use  the following convenient notation: by $\mathcal{O}(\varepsilon)$ we mean a positive function defined on $\Gamma$  such that $\displaystyle \lim_{\varepsilon \rightarrow 0} \mathcal{O}(\nu(\varepsilon)) \ = \ 0$.

\subsection{The axioms}

The first axiom is  a preparation for the next axioms. That is why we 
counted it as axiom 0.

\begin{enumerate}
\item[{\bf A0.}] The dilatations $$ \delta_{\varepsilon}^{x}: U(x) 
\rightarrow V_{\varepsilon}(x)$$ are defined for any 
$\displaystyle \varepsilon \in \Gamma, \nu(\varepsilon)\leq 1$. The sets 
$\displaystyle U(x), V_{\varepsilon}(x)$ are open neighbourhoods of $x$.  
All dilatations are homeomorphisms (invertible, continuous, with 
continuous inverse). 

We suppose  that there is a number  $1<A$ such that for any $x \in X$ we have 
$$\bar{B}_{d}(x,A) \subset U(x)  \ .$$
 We suppose that for all $\varepsilon \in \Gamma$, $\nu(\varepsilon) \in 
(0,1)$, we have 
$$ B_{d}(x,\nu(\varepsilon)) \subset \delta_{\varepsilon}^{x} B_{d}(x,A) 
\subset V_{\varepsilon}(x) \subset U(x) \ .$$

 There is a number $B \in (1,A]$ such that  for 
 any $\varepsilon \in \Gamma$ with $\nu(\varepsilon) \in (1,+\infty)$ the 
 associated dilatation  
$$\delta^{x}_{\varepsilon} : W_{\varepsilon}(x) \rightarrow B_{d}(x,B) \ , $$
is injective, invertible on the image. We shall suppose that 
$\displaystyle  W_{\varepsilon}(x) \in \mathcal{V}(x)$, that   
$\displaystyle V_{\varepsilon^{-1}}(x) \subset W_{\varepsilon}(x) $
and that for all $\displaystyle \varepsilon \in \Gamma_{1}$ and 
$\displaystyle u \in U(x)$ we have
$$\delta_{\varepsilon^{-1}}^{x} \ \delta^{x}_{\varepsilon} u \ = \ u \ .$$
\end{enumerate}

We have therefore  the following string of inclusions, for any $\varepsilon \in \Gamma$, $\nu(\varepsilon) \leq 1$, and any $x \in X$:
$$ B_{d}(x,\nu(\varepsilon)) \subset \delta^{x}_{\varepsilon}  B_{d}(x, A) 
\subset V_{\varepsilon}(x) \subset 
W_{\varepsilon^{-1}}(x) \subset \delta_{\varepsilon}^{x}  B_{d}(x, B) \quad . $$

A further technical condition on the sets  $\displaystyle V_{\varepsilon}(x)$ and $\displaystyle W_{\varepsilon}(x)$  will be given just before the axiom A4. (This condition will be counted as part of 
axiom A0.)

\begin{enumerate}
\item[{\bf A1.}]  We  have 
$\displaystyle  \delta^{x}_{\varepsilon} x = x $ for any point $x$. We also have $\displaystyle \delta^{x}_{1} = id$ for any $x \in X$.

Let us define the topological space
$$ dom \, \delta = \left\{ (\varepsilon, x, y) \in \Gamma \times X \times X 
\mbox{ : } \quad \mbox{ if } \nu(\varepsilon) \leq 1 \mbox{ then } y \in U(x) \,
\, , 
\right.$$ 
$$\left. \mbox{  else } y \in W_{\varepsilon}(x) \right\} $$ 
with the topology inherited from the product topology on 
$\Gamma \times X \times X$. Consider also $\displaystyle Cl(dom \, \delta)$, 
the closure of $dom \, \delta$ in $\displaystyle \bar{\Gamma} \times X \times X$ with product topology. 
The function $\displaystyle \delta : dom \, \delta \rightarrow  X$ defined by 
$\displaystyle \delta (\varepsilon,  x, y)  = \delta^{x}_{\varepsilon} y$ is continuous. Moreover, it can be continuously extended to $\displaystyle Cl(dom \, \delta)$ and we have 
$$\lim_{\varepsilon\rightarrow 0} \delta_{\varepsilon}^{x} y \, = \, x \quad . $$

\item[{\bf A2.}] For any  $x, \in K$, $\displaystyle \varepsilon, \mu \in \Gamma_{1}$ and $\displaystyle u \in 
\bar{B}_{d}(x,A)$   we have: 
$$ \delta_{\varepsilon}^{x} \delta_{\mu}^{x} u  = \delta_{\varepsilon \mu}^{x} u  \ .$$

\item[{\bf A3.}]  For any $x$ there is a  function $\displaystyle (u,v) \mapsto d^{x}(u,v)$, defined for any $u,v$ in the closed ball (in distance d) $\displaystyle 
\bar{B}(x,A)$, such that 
$$\lim_{\varepsilon \rightarrow 0} \quad \sup  \left\{  \mid \frac{1}{\varepsilon} d(\delta^{x}_{\varepsilon} u, \delta^{x}_{\varepsilon} v) \ - \ d^{x}(u,v) \mid \mbox{ :  } u,v \in \bar{B}_{d}(x,A)\right\} \ =  \ 0$$
uniformly with respect to $x$ in compact set. 

\end{enumerate}

\begin{remark}
The "distance" $d^{x}$ can be degenerated: there might exist  
$\displaystyle v,w \in U(x)$ such that $\displaystyle d^{x}(v,w) = 0$. 
\label{imprk}
\end{remark}

For  the following axiom to make sense we impose a technical condition on the co-domains $\displaystyle V_{\varepsilon}(x)$: for any compact set $K \subset X$ there are $R=R(K) > 0$ and 
$\displaystyle \varepsilon_{0}= \varepsilon(K) \in (0,1)$  such that  
for all $\displaystyle u,v \in \bar{B}_{d}(x,R)$ and all $\displaystyle \varepsilon \in \Gamma$, $\displaystyle  \nu(\varepsilon) \in (0,\varepsilon_{0})$,  we have 
$$\delta_{\varepsilon}^{x} v \in W_{\varepsilon^{-1}}( \delta^{x}_{\varepsilon}u) \ .$$

With this assumption the following notation makes sense:
$$\Delta^{x}_{\varepsilon}(u,v) = \delta_{\varepsilon^{-1}}^{\delta^{x}_{\varepsilon} u} \delta^{x}_{\varepsilon} v . $$
The next axiom can now be stated: 
\begin{enumerate}
\item[{\bf A4.}] We have the limit 
$$\lim_{\varepsilon \rightarrow 0}  \Delta^{x}_{\varepsilon}(u,v) =  \Delta^{x}(u, v)  $$
uniformly with respect to $x, u, v$ in compact set. 
\end{enumerate}

\begin{definition}
A triple $(X,d,\delta)$ which satisfies A0, A1, A2, A3, but $\displaystyle d^{x}$ is degenerate for some 
$x\in X$, is called degenerate dilatation structure. 

If the triple $(X,d,\delta)$ satisfies A0, A1, A2, A3 and 
 $\displaystyle d^{x}$ is non-degenerate for any $x\in X$, then we call it  a 
 dilatation structure. 

 If a  dilatation structure satisfies A4 then we call it strong dilatation 
 structure. 
 \label{defweakstrong}
\end{definition}

\section{Dilatation structures on the boundary of the dyadic tree}

Dilatation structures on the boundary of the dyadic tree will have a simpler 
form than general, mainly because the distance is ultrametric.

We shall take the group $\Gamma$ to be the set of integer powers of $2$, seen as a subset of 
dyadic numbers. Thus for any $p \in \mathbb{Z}$ the element $\displaystyle 2^{p} \in \mathbb{Q}_{2}$ 
belongs to $\Gamma$. The operation is the multiplication of dyadic numbers and the morphism 
$\nu : \Gamma \rightarrow (0,+\infty)$ is defined by 
$$\nu(2^{p}) = d(0, 2^{p}) = \frac{1}{2^{p}} \in (0,+\infty) \quad  . $$

\paragraph{Axiom A0.} This axiom  states that for any $p\in \mathbb{N}$ and any $x \in X^{\omega}$ the dilatation 
$$\delta^{x}_{2^{p}} : U(x) \rightarrow V_{2^{p}}(x) $$
is a homeomorphism, the sets $U(x)$ and $\displaystyle V_{2^{p}}(x) $ are open and 
there is $A>1$ such that the ball centered in $x$ and radius $A$ is contained in $U(x)$. But this means 
that $\displaystyle U(x) = X^{\omega}$, because $\displaystyle X^{\omega} = B(x,1)$. 

Further, for any $p \in \mathbb{N}$ we have the inclusions: 
\begin{equation}
B(x, \frac{1}{2^{p}}) \subset \delta^{x}_{2^{p}} X^{\omega} \subset V_{2^{p}}(x) \quad . 
\label{a01}
\end{equation}

For any $\displaystyle p \in \mathbb{N}^{*}$ the associated dilatation  
$\displaystyle \delta^{x}_{2^{-p}} : W_{2^{-p}}(x) \rightarrow B(x,B) = 
X^{\omega}$ , 
is injective, invertible on the image. We suppose that 
$\displaystyle W_{2^{-p}}(x)$ is open, that 
\begin{equation}
V_{2^{p}}(x) \subset W_{2^{-p}}(x) 
\label{a02}
\end{equation}
and that for all $\displaystyle p \in \mathbb{N}^{*}$ and $\displaystyle u \in X^{\omega}$ we have
$\displaystyle \delta_{2^{-p}}^{x} \ \delta^{x}_{2^{p}} u \ = \ u$ . 
We leave aside for the moment the interpretation of the technical condition before axiom A4.
 
\paragraph{Axioms A1 and A2.}  Nothing simplifies.

\paragraph{Axiom A3.} Because $d$ is an ultrametric distance and $\displaystyle X^{\omega}$ is compact, this axiom has very strong consequences, for a non degenerate dilatation structure. 

In this case the axiom A3 states that there is a non degenerate distance function  
$\displaystyle d^{x}$ on 
$\displaystyle X^{\omega}$ such that we have the limit 
\begin{equation}
\lim_{p \rightarrow \infty} 2^{p} d( \delta^{x}_{2^{p}} u, \delta^{x}_{2^{p}} v) = d^{x}(u,v) 
\label{a30}
\end{equation}
uniformly with respect to $\displaystyle x,u,v \in X^{\omega}$. 

We continue further with first properties of  dilatation structures. 

\begin{lemma}
There exists $\displaystyle p_{0} \in \mathbb{N}$ such that for any $\displaystyle x, u, 
v \in X^{\omega}$ and for any $\displaystyle p \in \mathbb{N}$, $\displaystyle p \geq  p_{0}$, we have 
$$ 2^{p} d( \delta^{x}_{2^{p}} u, \delta^{x}_{2^{p}} v) = d^{x}(u,v) \quad  . $$
\label{l1}
\end{lemma}

\begin{proof}
From the limit (\ref{a30}) and the non degeneracy of the distances $\displaystyle d^{x}$ we deduce that 
$$\lim_{p \rightarrow \infty} \log_{2} \left( 2^{p} d( \delta^{x}_{2^{P}} u, \delta^{x}_{2^{P}} v)\right)  =  
\log_{2} d^{x}(u,v) \quad , $$
uniformly with respect to $\displaystyle x,u,v \in X^{\omega}$, $u \not = v$. The right hand side term is finite and the sequence from the limit at the left hand side is included in $\mathbb{Z}$. Use this and 
the uniformity of the convergence to get the desired result. 
\end{proof}

 
 In the sequel $\displaystyle p_{0}$ is the smallest natural number satisfying lemma \ref{l1}. 
 
 \begin{lemma}
For any $\displaystyle x  \in X^{\omega}$ and for any 
$\displaystyle p \in \mathbb{N}$, $\displaystyle p \geq  p_{0}$, we have 
$\displaystyle   \delta^{x}_{2^{p}} X^{\omega} = [x]_{p} X^{\omega}$. 
Otherwise stated, for any $\displaystyle x, y \in X^{\omega}$, any $\displaystyle q  \in X^{*}$, 
$\displaystyle \mid q \mid \geq p_{0}$ there exists $\displaystyle w  \in X^{\omega}$ such that 
 $$  \delta^{qx}_{2^{\mid q \mid}} w = qy $$
 and for any $\displaystyle z \in X^{\omega}$ there is  
 $\displaystyle y \in X^{\omega}$ such that 
 $\displaystyle   \delta^{qx}_{2^{\mid q \mid}} z = qy$ . 
 Moreover, for any $\displaystyle x  \in X^{\omega}$ and for any $\displaystyle p \in \mathbb{N}$, $\displaystyle p \geq  p_{0}$ the inclusions from (\ref{a01}), (\ref{a02}) are equalities. 
\label{l2}
\end{lemma}

\begin{proof}
From the last inclusion in  (\ref{a01})  we get that for any 
$\displaystyle x, y \in X^{\omega}$, any $\displaystyle q  \in X^{*}$, 
$\displaystyle \mid q \mid \geq p_{0}$ there exists $\displaystyle w  \in X^{\omega}$ such that 
 $\displaystyle  \delta^{qx}_{2^{\mid q \mid}} w = qy$ . 
 For the second part of the conclusion we use lemma \ref{l1} and axiom A1. From there we see that 
 for any $\displaystyle p \geq p_{0}$ we have 
 $$2^{p} d( \delta^{x}_{2^{p}} x, \delta^{x}_{2^{p}} u) =  2^{p} d( x, \delta^{x}_{2^{p}} u) = d^{x}(x,u) 
 \leq 1\quad  . $$
 Therefore $\displaystyle 2^{p} d( x, \delta^{x}_{2^{p}} u) \leq 1$, which is equivalent with the second part of the lemma. 
 
 Finally, the last part of the lemma has a similar proof, only that we have to use also the last part of axiom A0. 
 \end{proof}
 

The technical condition before the axiom A4 turns out to be trivial. Indeed, from lemma \ref{l2} it follows that for any $\displaystyle p \geq p_{0}$, $p \in \mathbb{N}$, and any $\displaystyle x, u, v \in 
X^{\omega}$  we have $\displaystyle \delta^{x}_{2^{p}} u = [x]_{p} w$, $\displaystyle w \in X^{\omega}$.  It follows that  
$$\delta^{x}_{2^{p}} v \in [x]_{p} X^{\omega} = W_{2^{-p}} (x) = W_{2^{-p}} (  \delta^{x}_{2^{p}} u) \quad . $$

 \begin{lemma}
For any $\displaystyle x, u, 
v \in X^{\omega}$  such that $\displaystyle 2^{p_{0}} d(x,u) \leq 1$,  
$\displaystyle 2^{p_{0}} d(x,v) \leq 1$ we have 
$\displaystyle d^{x}(u,v) = d(u,v)$ . 
Moreover, under the same hypothesis,  for any $\displaystyle p \in \mathbb{N}$ we have 
$$ 2^{p} d( \delta^{x}_{2^{p}} u, \delta^{x}_{2^{p}} v) = d(u,v) \quad  . $$
\label{l3}
\end{lemma}

\begin{proof}
By lemma \ref{l1}, lemma \ref{l2} and axiom A2. Indeed, from lemma \ref{l1} and axiom A2,  for any $\displaystyle p \in \mathbb{N}$ and any 
$\displaystyle x, u', v' \in X^{\omega}$ we have 
$$d^{x}(u',v') = 2^{p_{0}+p} d( \delta^{x}_{2^{p+p_{0}}} u' ,  \delta^{x}_{2^{p+p_{0}}} v' ) = $$ 
$$ = 2^{p} \,   2^{p_{0}} d( \delta^{x}_{2^{p_{0}}}   \delta^{x}_{2^{p}}u', \delta^{x}_{2^{p_{0}}}   \delta^{x}_{2^{p}}v') = 2^{p} d^{x}( \delta^{x}_{2^{p}}u' , \delta^{x}_{2^{p}}v') \quad . $$
This is just the cone property for $\displaystyle d^{x}$. From here we deduce 
that for any $\displaystyle p \in \mathbb{Z}$ we have 
$\displaystyle d^{x}(u', v') = 2^{p} d^{x}(\delta^{x}_{2^{p}}u' ,
\delta^{x}_{2^{p}}v')$ . 
If $\displaystyle 2^{p_{0}} d(x,u) \leq 1$,  $\displaystyle 2^{p_{0}} d(x,v) 
\leq 1$  then write $x = qx'$, $\displaystyle \mid q \mid = p_{0}$, and use 
lemma \ref{l2} to get the existence of $\displaystyle u', v' \in X^{\omega}$ 
such that $\displaystyle \delta^{x}_{2^{p_{0}}} u' = u$ , $\displaystyle 
\delta^{x}_{2^{p_{0}}} v' = v$ . 
Therefore, by lemma \ref{l1}, we have 
$$ d(u,v) = 2^{-p_{0}} d^{x}( u', v') = d^{x}( \delta^{x}_{2^{-p_{0}}} u' ,   \delta^{x}_{2^{-p_{0}}} v' ) = 
d^{x}(u,v) \quad . $$
The first part of the lemma is proven. For the proof of the second part write again 
$$2^{p} d( \delta^{x}_{2^{p}} u, \delta^{x}_{2^{p}} v) = 2^{p} d^{x}( \delta^{x}_{2^{p}} u, \delta^{x}_{2^{p}} v) = d^{x}(u,v) = d(u,v) 
$$
which finishes the proof. 
\end{proof}

The space $\displaystyle X^{\omega}$ decomposes into a disjoint union of 
$\displaystyle 2^{p_{0}}$ balls which are isometric. There is no connection 
between the  dilatation structures on these balls, therefore we shall 
suppose further that $\displaystyle p_{0} = 0$.

Our purpose is to find the  general form  of a 
dilatation structure   on $\displaystyle X^{\omega}$,  with  
$\displaystyle p_{0} = 0$. 

\begin{definition}
A function $\displaystyle W: \mathbb{N}^{*} \times X^{\omega} \rightarrow 
Isom(X^{\omega})$  is smooth if for any $\varepsilon > 0$ there exists $\mu(\varepsilon) > 0$ such that 
for any $\displaystyle x, x' \in X^{\omega}$ such that $d(x,x')< \mu(\varepsilon)$ and for any $\displaystyle y \in X^{\omega}$  we have 
$$ \frac{1}{2^{k}} \, d( W^{x}_{k} (y) , W^{x'}_{k} (y) ) \leq \varepsilon \quad , $$
for an  $k$ such that  $\displaystyle d(x,x') <   1 / 2^{k} $.  
\label{defwsmooth}
\end{definition}

\begin{theorem}
Let $\displaystyle (X^{\omega}, d, \delta)$ be a  dilatation structure on $\displaystyle (X^{\omega}, d)$, where $d$ is the standard distance on $\displaystyle X^{\omega}$, such that $\displaystyle p_{0} = 0$. Then there exists a  smooth (according to definition \ref{defwsmooth}) function 
$$W: \mathbb{N}^{*} \times X^{\omega} \rightarrow Isom(X^{\omega}) \quad , \quad W(n,x) = W^{x}_{n} $$
such that  for any 
$\displaystyle q \in X^{*}$, $\alpha \in X$, $x, y \in X^{\omega}$ we have 
\begin{equation}
\delta_{2}^{q \alpha x} q \bar{\alpha} y = q \alpha \bar{x_{1}} W^{q \alpha x}_{\mid q \mid + 1} 
(y) \quad  . 
\label{eqtstruc}
\end{equation}
Conversely, to any smooth function  $\displaystyle W: \mathbb{N}^{*} \times X^{\omega} \rightarrow 
Isom(X^{\omega})$ is associated a  dilatation structure  $\displaystyle (X^{\omega}, d, \delta)$, 
with $\displaystyle p_{0} = 0$, induced by  functions $\displaystyle \delta_{2}^{x}$, defined by 
$\displaystyle \delta_{2}^{x} x = x$ and otherwise by relation (\ref{eqtstruc}).
\label{tstruc}
\end{theorem}
 
\begin{proof}
Let $\displaystyle (X^{\omega}, d, \delta)$ be a   dilatation structure on 
$\displaystyle (X^{\omega}, d)$, such that $\displaystyle p_{0} = 0$. 
Any two different elements of $\displaystyle X^{\omega}$ can be written in the 
form $q \alpha x$ and $q \bar{\alpha} y$, with $\displaystyle q \in X^{*}$, 
$\alpha \in X$, $x, y \in X^{\omega}$. We also have 
$\displaystyle d(q \alpha x , q \bar{\alpha} y) = 2^{- \mid q \mid}$ . 
From the following computation (using $\displaystyle p_{0} = 0$ and axiom A1): 
$$2^{-\mid q \mid -1} = \frac{1}{2} d(q \alpha x , q \bar{\alpha} y) = d( q 
\alpha x ,  \delta_{2}^{q \alpha x} q \bar{\alpha} y) \quad , $$
we find that there exists $\displaystyle w^{q \alpha x}_{\mid q \mid + 1}(y)
 \in X^{\omega}$ such that 
$\displaystyle \delta_{2}^{q \alpha x} q \bar{\alpha} y = q \alpha 
w^{q \alpha x}_{\mid q \mid + 1}(y)$ . 
Further on, we compute: 
$$ \frac{1}{2} d(q \bar{\alpha} x , q \bar{\alpha} y) = d( \delta_{2}^{q \alpha x} q \bar{\alpha} x , \delta_{2}^{q \alpha x} q \bar{\alpha} y) = d( q \alpha w^{q \alpha x}_{\mid q \mid + 1}(x) , q \alpha w^{q \alpha x}_{\mid q \mid + 1}(y)) \quad . $$
From this equality we find that 
$\displaystyle 1 > \frac{1}{2} d(x,y) = d( w^{q \alpha x}_{\mid q \mid + 1}(x) 
, w^{q \alpha x}_{\mid q \mid + 1}(y))$ , 
which means that the first letter of the word $\displaystyle w^{q \alpha x}_{\mid q \mid + 1}(y)$ 
does not depend on $y$, and is equal to the first letter of the word $\displaystyle w^{q \alpha x}_{\mid q \mid + 1}(x)$. Let us denote this letter by $\beta$ (which depends only on $q$, $\alpha$, $x$). Therefore we may write: 
$$w^{q \alpha x}_{\mid q \mid + 1}(y) = \beta W^{q \alpha x}_{\mid q \mid + 1}(y) \quad , $$
where the properties of the function $\displaystyle y \mapsto W^{q \alpha x}_{\mid q \mid + 1}(y)$ remain to be determined later. 

We go back to the first  computation in this proof: 
 $$2^{-\mid q \mid -1} = d( q \alpha x ,  \delta_{2}^{q \alpha x} q \bar{\alpha} y) = d( q \alpha x , q \alpha 
 \beta W^{q \alpha x}_{\mid q \mid + 1}(y)) \quad . $$
This shows that $\displaystyle \bar{\beta}$ is the first letter of the word $x$.  We proved the relation 
(\ref{eqtstruc}), excepting the fact that the function $\displaystyle y \mapsto W^{q \alpha x}_{\mid q \mid + 1}(y)$ is an isometry. But this is true. Indeed, for any $\displaystyle u,v \in X^{\omega}$ we have 
$$\frac{1}{2} d(q \bar{\alpha} u , q \bar{\alpha} v) = d( \delta_{2}^{q \alpha x} q \bar{\alpha} u , 
\delta_{2}^{q \alpha x} q \bar{\alpha} v) = d( q \alpha \bar{x_{1}} W^{q \alpha x}_{\mid q \mid + 1}(x) , q \alpha \bar{x_{1}} W^{q \alpha x}_{\mid q \mid + 1}(y)) \quad . $$
This proves the isometry property. 

The dilatations of coefficient $2$ induce all dilatations (by axiom A2). 
In order to satisfy the continuity 
assumptions from axiom A1, the function 
$\displaystyle W: \mathbb{N}^{*} \times X^{\omega} \rightarrow 
Isom(X^{\omega})$ has to be smooth in the sense of definition 
\ref{defwsmooth}.  Indeed, axiom A1 is equivalent to the fact that  
$\displaystyle \delta^{x'}_{2}(y')$ converges uniformly to $
\displaystyle \delta^{x}_{2}(y)$, as $d(x,x'), d(y,y')$ go to zero. 
There are two cases to study. 

Case 1: $d(x,x') \leq d(x,y)$, $d(y,y') \leq d(x,y)$. It means that $x = q \alpha q' \beta X$, $y = q \bar{\alpha} q" \gamma Y$, $x' = q \alpha q' \bar{\beta} X'$, $y' = q \bar{\alpha} q" \bar{\gamma} Y'$, 
with $d(x,y) = 1 / 2^{k}$, $k = \mid q \mid$ .

Suppose that $q' \not = \emptyset$. We compute then: 
$\displaystyle  \delta_{2}^{x}(y) = q \alpha \bar{q'_{1}} W_{k+1}^{x}( q" 
\gamma Y)$ , $\displaystyle  \delta_{2}^{x'}(y') = 
q \alpha \bar{q'_{1}} W_{k+1}^{x'}(q" \bar{\gamma} Y')$. 
All the functions denoted by a capitalized "W" are isometries, therefore we get the estimation: 
$$d(  \delta_{2}^{x}(y) ,  \delta_{2}^{x'}(y')) = \frac{1}{2^{k+2}} \, d(  W_{k+1}^{x}( q" \gamma Y) , W_{k+1}^{x'}(q" \bar{\gamma} Y')) \leq $$ 
$$ \leq \frac{1}{2^{k+2}} \, d( q" \gamma Y , q" \bar{\gamma} Y') + 
\, \frac{1}{2^{k+2}} \,  d( W_{k+1}^{x}(q" \gamma Y) , W_{k+1}^{x'}(q" \gamma Y)) = $$
$$ =  \frac{1}{2} d(y,y') + 
\, \frac{1}{2^{k+2}} \, d( W_{k+1}^{x}(q" \gamma Y) , W_{k+1}^{x'}(q" \gamma Y)) \quad . $$
We see that if $W$ is smooth in the sense of definition \ref{defwsmooth} then 
the structure $\delta$ satisfies the uniform continuity assumptions for this 
case. Conversely, if $\delta$ satisfies A1 then 
$W$ has to be smooth. 

If $q' = \emptyset$ then a similar computation leads to the same conclusion. 

Case 2: $d(x,x') > d(x,y) > d(y,y')$. It means that $x = q \alpha q' \beta X$, $x' = q \bar{\alpha} X'$, 
$y = q \alpha q' \bar{\beta} q" \bar{\gamma} Y$, $y' =  q \alpha q' \bar{\beta} q" \gamma Y'$, with 
$d(x,x') = 1 \ 2^{k}$, $k = \mid q \mid$ .  

We compute then: 
$\displaystyle  \delta_{2}^{x}(y) = q \alpha q ' \beta \bar{X_{1}} 
W_{k+2+ \mid q' \mid}^{x}( q" \bar{\gamma} Y)$, 
$$\delta_{2}^{x'}(y') = 
q \bar{\alpha} \bar{X'_{1}} W_{k+1}^{x'}(q' \bar{\beta} q" \gamma Y') \leq  
 \frac{1}{2^{k}} = d(x, x') \quad . $$
Therefore in his case the continuity is satisfied, without any supplementary constraints on the function $W$. 

The first part of the theorem is proven. 

For the proof of the second part of the theorem we start from the function 
$\displaystyle W: \mathbb{N}^{*} \times X^{\omega} \rightarrow Isom(X^{\omega})$. It is sufficient to prove  for any $\displaystyle 
x, y, z \in X^{\omega}$ the equality 
$$\frac{1}{2} d(y,z) = d(\delta_{2}^{x} y, \delta_{2}^{x} z) \quad .$$
Indeed, then we can construct the all dilatations from the dilatations of 
coefficient $2$ (thus we satisfy A2). All axioms, excepting A1, are satisfied. 
But A1 is equivalent with the smoothness of the function 
$W$, as we proved earlier.

Let us prove now the  before mentioned equality. 
If $y = z$ there is nothing to prove.  Suppose that $y \not = z$. The distance $d$ is ultrametric, therefore the proof splits in two cases. 

Case 1: $d(x,y) = d(x,z) > d(y,z)$. This is equivalent to $x = q \bar{\alpha} x'$, $y = q \alpha q' \beta y'$, 
$z = q \alpha q' \bar{\beta} z'$, with $\displaystyle q, q' \in X^{*}$, $\alpha, \beta \in X$, 
$\displaystyle x', y', z' \in X^{\omega}$. We compute: 
$$d( \delta_{2}^{x} y, \delta_{2}^{x} z) = d(\delta_{2}^{q \bar{\alpha} x'} q \alpha q' \beta y', \delta_{2}^{q \bar{\alpha} x'} q \alpha q' \bar{\beta} z') = $$
$$= d(q \bar{\alpha} \bar{x'_{1}} W^{x}_{\mid q \mid +1} ( q' \beta y') ,  q \bar{\alpha} \bar{x'_{1}} W^{x}_{\mid q \mid +1} ( q' \bar{\beta} z')) = 2^{-\mid q \mid - 1} d(  W^{x}_{\mid q \mid +1} ( q' \beta y'), 
 W^{x}_{\mid q \mid +1} ( q' \bar{\beta} z')) = $$
 $$ = 2^{-\mid q \mid - 1}  d(q' \beta y', q' \bar{\beta} z')) = \frac{1}{2} d(q \alpha q' \beta y' , q \alpha q' \bar{\beta} z') = \frac{1}{2} d(y,z) \quad  . $$

 Case 2: $d(x,y) = d(y,z) > d(x,z)$. If $x = z$ then we write $x = q \alpha u$, $y = q \bar{\alpha} v$ and we have 
 $$  d(\delta_{2}^{x} y, \delta_{2}^{x} z) = d( q \alpha \bar{u_{1}} W^{x}_{\mid q \mid +1} (v), q \alpha u) 
 = 2^{-\mid q \mid + 1} =  \frac{1}{2} d(y,z) \quad . $$
If $x \not = z$ then we can write $z = q \bar{\alpha} z'$, $y = q \alpha q' \beta y'$, 
$x = q \alpha q' \bar{\beta} x'$, with $\displaystyle q, q' \in X^{*}$, $\alpha, \beta \in X$, 
$\displaystyle x', y', z' \in X^{\omega}$. We compute: 
$$d( \delta_{2}^{x} y, \delta_{2}^{x} z) = d(\delta_{2}^{q \alpha q' \bar{\beta} x'}  q \alpha q' \beta y' , 
\delta_{2}^{q \alpha q' \bar{\beta} x'} q \bar{\alpha} z') = $$
$$ = d(q \alpha q' \bar{\beta} \bar{x'_{1}} W^{x}_{\mid q \mid + \mid q' \mid +2}(y'), q \alpha \gamma 
W^{x}_{\mid q \mid + 1} (z')) \quad ,  $$
with $\gamma \in X$, $\displaystyle \bar{\gamma} = q'_{1}$ if $q' \not = \emptyset$, otherwise 
$\gamma = \beta$. In both situations we have 
$\displaystyle d( \delta_{2}^{x} y, \delta_{2}^{x} z) =  
2^{- \mid q \mid - 1} = \frac{1}{2} d(y,z)$ . 
The proof is done. 
\end{proof}


\subsection{Self-similar  dilatation structures}

Let $\displaystyle (X^{\omega}, d, \delta)$ be a  dilatation structure. There are induced 
dilatations structures on $\displaystyle 0X^{\omega}$ and  $\displaystyle 1X^{\omega}$. 

\begin{definition}
For any $\alpha \in X$ and $\displaystyle x, y \in X^{\omega}$ we define  $\displaystyle \delta_{2}^{\alpha, x} y$ by the relation 
$$\delta_{2}^{\alpha x} \alpha y = \alpha \, \delta_{2}^{\alpha, x} y \quad . $$
\end{definition}

The following proposition has a straightforward proof, therefore we skip it. 

\begin{proposition}
If $\displaystyle (X^{\omega}, d, \delta)$ is a   dilatation structure and $\alpha \in X$ then 
$\displaystyle (X^{\omega}, d, \delta^{\alpha})$ is a   dilatation structure. 

If $\displaystyle (X^{\omega}, d, \delta')$ and $\displaystyle (X^{\omega}, d, \delta")$ are   dilatation structures then $\displaystyle (X^{\omega}, d, \delta)$ is a   dilatation structure, where 
$\delta$ is uniquely defined by $\displaystyle \delta^{0} = \delta'$, $\displaystyle \delta^{1} = \delta"$. 
\end{proposition}

\begin{definition}
 A   dilatation structure $\displaystyle (X^{\omega}, d, \delta)$ is self-similar if for any $\alpha \in X$ and $\displaystyle x,y \in X^{\omega}$ we have 
 $$\delta_{2}^{\alpha x} \, \alpha y = \alpha \, \delta_{2}^{x} y \quad . $$
 \label{defself}
\end{definition}

 Self-similarity is thus related to linearity. Indeed, let us compare 
 self-similarity with the following definition of linearity. 
 
 \begin{definition}
 For a given dilatation structure $\displaystyle (X^{\omega}, d, \delta)$, a 
 continuous transformation $\displaystyle A: X^{\omega} \rightarrow X^{\omega}$ is linear (with respect 
 to the dilatation structure) if for any $\displaystyle x , y \in X^{\omega}$ we
 have 
 $$A \, \delta^{x}_{2} y \, = \, \delta^{A x}_{2} A y$$
 \label{defline}
 \end{definition}
 
 The previous definition provides a true generalization of linearity for 
 dilatation structures. This can be seen by comparison with the characterisation
 of linear (in fact affine) transformations in vector spaces from the 
  proposition \ref{1proplin}.  
  
 The definition of self-similarity \ref{defself} is related to linearity in 
 the sense of definition \ref{defline}. To see this, let us consider the
 functions $\displaystyle \hat{\alpha} : X^{\omega} \rightarrow X^{\omega}$, 
 $\displaystyle \hat{\alpha} x \, = \, \alpha x$, for $\alpha \in X$. With this 
 notations, the definition \ref{defself} simply states that a dilatation
 structure is self-similar if these two functions, $\displaystyle \hat{0}$ and 
 $\displaystyle \hat{1}$, are linear in the sense of definition \ref{defline}.

 The description of self-similar dilatation structures on the boundary of the
 dyadic tree is given in the next theorem. 

\begin{theorem}
Let $\displaystyle (X^{\omega}, d, \delta)$ be a self-similar   dilatation structure and 
$\displaystyle W: \mathbb{N}^{*} \times X^{\omega} \rightarrow Isom(X^{\omega})$ the function 
associated to it, according to theorem \ref{tstruc}.  Then there exists a 
function $W: X^{\omega} \rightarrow Isom(X^{\omega})$ such that: 
\begin{enumerate}
\item[(a)] for any $\displaystyle q \in X^{*}$ and any $\displaystyle x \in X^{\omega}$ we have 
$\displaystyle W_{\mid q \mid + 1}^{q x} = W^{x}$ , 
\item[(b)] there exists $C>0$ such that for any $\displaystyle x, x', y \in X^{\omega}$ and for any 
$\lambda > 0$, if $d(x,x') \leq \lambda$ then 
$\displaystyle d(W^{x}(y) , W^{x'}(y)) \leq C \lambda$ . 
\end{enumerate}
\label{th2}
\end{theorem}

\begin{proof}
We define $\displaystyle W^{x} = W^{x}_{1}$ for any $\displaystyle x \in X^{\omega}$ .  
We want to prove that this function satisfies (a), (b). 

(a) Let $\displaystyle \beta \in X$ and any $\displaystyle x, y \in X^{\omega}$, 
$x = q \alpha u$, $y = q \bar{\alpha} v$. By self-similarity we obtain: 
$\displaystyle \beta q \alpha \bar{u_{1}} W_{\mid q \mid + 2}^{\beta x} (v) = 
\delta_{2}^{\beta x} \beta y =  \beta \delta_{2}^{x} y =  
\beta q \alpha \bar{u_{1}} W_{\mid q \mid + 1}^{x} (v)$ . 
We proved that 
$$ W_{\mid q \mid + 2}^{\beta x} (v) = W_{\mid q \mid + 1}^{x} (v)  $$
for any $\displaystyle x, v \in X^{\omega}$ and $\beta \in X$ This implies (a).  

(b) This is a consequence of smoothness,  in the sense of definition \ref{defwsmooth}, of the function 
$\displaystyle W: \mathbb{N}^{*} \times X^{\omega} \rightarrow Isom(X^{\omega})$.  Indeed, 
$\displaystyle (X^{\omega}, d, \delta)$ is a    dilatation structure, therefore by theorem \ref{tstruc} 
the  previous mentioned function is smooth. 

By (a) the smoothness condition becomes: for any $\varepsilon > 0$ there is $\mu(\varepsilon) > 0$ such that for any $\displaystyle y \in X^{\omega}$, any $k \in \mathbb{N}$ and  any $x, x' \in X^{\omega}$,  
if $d(x,x') \leq 2^{k} \mu(\varepsilon)$ then 
$$d( W^{x}(y), W^{x'} (y)) \leq 2^{k} \varepsilon \quad . $$
Define then the modulus  of continuity: for any $\varepsilon > 0$ let $\bar{\mu}(\varepsilon)$ be given by 
$$\bar{\mu}(\varepsilon) = \sup \left\{ \mu \, \mbox{ : } \forall x, x', y \in X^{\omega} \, \, d(x,x')\leq \mu \Longrightarrow d( W^{x}(y), W^{x'} (y)) \leq \varepsilon \right\} \quad .$$
We see that the modulus of continuity $\bar{\mu}$ has the property 
$$\bar{\mu}(2^{k} \varepsilon)  = 2^{k} \bar{\mu}(\varepsilon) $$
for any $k \in \mathbb{N}$. Therefore there exists $C> 0$ such that 
$\displaystyle \bar{\mu}( \varepsilon) = C^{-1} \varepsilon$ for any $\displaystyle \varepsilon = 1/2^{p}$, $p \in \mathbb{N}$. The point (b) follows immediately. 
\end{proof}


\end{document}